\documentclass{article}
\usepackage{amsmath, amsthm}
\usepackage{amssymb}
\newtheoremstyle{theorem}
  {10pt}          
  {10pt}  
  {\sl}  
  {\parindent}     
  {\bf}  
  {. }    
  { }    
  {}     
\newtheorem{thm}{Theorem}[section]

\newtheorem{lem}[thm]{Lemma}
\newtheorem{prop}[thm]{Proposition}
\newtheorem{defn}{Definition}[section]

\newtheorem{rem}{Remark}[section]
\begin{document}
\title{\large\bf Dual identities in fractional difference calculus within Riemann}
\author{\small \bf Thabet Abdeljawad $^a$'$^b$  \\ {\footnotesize $^a$ Department of
Mathematics, \c{C}ankaya University, 06530, Ankara, Turkey}\\ {\footnotesize $^b$ Department of Mathematics and Physical Sciences}\\
{\footnotesize Prince Sultan University, P. O. Box 66833, Riyadh 11586, Saudi Arabia}}
\date{}
\maketitle

{\footnotesize {\noindent\bf Abstract.} We Investigate two types of dual identities for Riemann  fractional sums and  differences. The first type relates nabla and delta type fractional sums and differences. The second type represented by the Q-operator relates left and right fractional sums and differences. These dual identities insist that in the definition of right fractional differences we have to use both the nabla and delta operators. The solution representation for higher order Riemann fractional difference equation is obtained as well.
 \\


\section{Introduction}

During the last two decades, due to its widespread applications in different fields of science and engineering, fractional calculus has attracted the attention of many researchers \cite{podlubny,Samko, Kilbas,agrawal1,scalas,tr}.

 Starting from the idea of  discretizing the Cauchy integral formula, Miller and Ross \cite{Miller} and Gray and Zhang \cite{Gray}  obtained discrete versions of left type fractional integrals and derivatives, called fractional sums and differences. Fifteen years later, several authors  started to deal with discrete fractional calculus \cite{Th Caputo,Ferd,Feri,Nabla,Atmodel,Nuno,TDbyparts,Gdelta,Gnabla,Gfound}, benefiting from the theory of time scales originated by Hilger in 1988 (see \cite{Adv}).

 In this article, we summarize some of the results mentioned in the above references and add more in  right type fractional differences and the higher order. Throughout the article we almost agree with the previously presented definitions  except for the definition of right fractional difference. We shall figure out that these definitions seem to be more convenient than the previously presented ones by proving some dual identities. These identities fall into two kinds. The first relates nabla type fractional differences and sums to delta ones. The second kind represented by the Q-operator relates left and right fractional sums and differences. This setting enables us to get identities resembling better the ordinary fractional case. Along with the previously mentioned points we are able to fit a reasonable nabla integration by parts formula which remains in accordance with the one obtained in \cite{TDbyparts} but different  from those obtained in \cite{Atmodel} and \cite{Nuno}. The obtained dual identities are also used to obtain a delta integration by parts formula from the nabla one. The solution representation for the  higher order Riemann fractional difference equation is obtained as well and thus the result in \cite{THFer} is generalized. We shall see that the higher order Riemann fractional difference initial value problem of non-integer order needs only one initial condition which is not the case for the higher order difference equation (i.e of positive integer order).

 The article is organized as follows: The remaining part of this section contains  summary to some of the basic notations and definitions in delta and nabla calculus. Section 2 contains the definitions in the frame of delta and nabla fractional sums and differences in the Riemann sense. Moreover,
 some essential lemmas about the commutativity  of the different fractional sum operators with the usual difference operators are established. These lemmas are vital to proceed in the next sections. The third section contains some dual identities relating nabla and delta fractional sums and differences in the left and right cases. Using these dual identities, power formulae  for nabla left and right fractional sums and a commutative law for nabla left and right sums are obtained. In Section 4 the integration by parts formula for nabla fractional sums and differences obtained in \cite{THFer} is used by the help of the dual identities to obtain delta integration by parts formulas.
Section 5 is devoted to higher order initial fractional   difference value problems into Riemann. Finally, in Section 6  the Q-operator is used to relate left and right fractional sums in the nabla and delta case and hence relates delta and nabla Riemann fractional differences. The Q-dual identities obtained in this section expose the validity of the definition of delta and nabla right fractional differences where the delta and nabla operators are used in each of the definitions of delta and nabla right fractional differences.

For a natural number $n$, the fractional polynomial is defined by,

 \begin{equation} \label{fp}
 t^{(n)}=\prod_{j=0}^{n-1} (t-j)=\frac{\Gamma(t+1)}{\Gamma(t+1-n)},
 \end{equation}
where $\Gamma$ denotes the special gamma function and the product is
zero when $t+1-j=0$ for some $j$. More generally, for arbitrary
$\alpha$, define
\begin{equation} \label{fpg}
t^{(\alpha)}=\frac{\Gamma(t+1)}{\Gamma(t+1-\alpha)},
\end{equation}
where the convention that division at pole yields zero.
Given that the forward and backward difference operators are defined
by
\begin{equation} \label{fb}
\Delta f(t)=f(t+1)-f(t)\texttt{,}~\nabla f(t)=f(t)-f(t-1)
\end{equation}
respectively, we define iteratively the operators
$\Delta^m=\Delta(\Delta^{m-1})$ and $\nabla^m=\nabla(\nabla^{m-1})$,
where $m$ is a natural number.

 Here are some properties of the  factorial function.

\begin{lem} \label{pfp} (\cite{Ferd})
Assume the following factorial functions are well defined.

(i) $ \Delta t^{(\alpha)}=\alpha  t^{(\alpha-1)}$.

(ii) $(t-\mu)t^{(\mu)}= t^{(\mu+1)}$, where $\mu \in \mathbb{R}$.

(iii) $\mu^{(\mu)}=\Gamma (\mu+1)$.

(iv) If $t\leq r$, then $t^{(\alpha)}\leq r^{(\alpha)}$ for any
$\alpha>r$.

(v) If $0<\alpha<1$, then $ t^{(\alpha\nu)}\geq
(t^{(\nu)})^\alpha$.

(vi) $t^{(\alpha+\beta)}= (t-\beta)^{(\alpha)} t^{(\beta)}$.
\end{lem}

Also, for our purposes we list down the following two properties, the proofs of which are straightforward.

\begin{equation} \label{ou1}
\nabla_s (s-t)^{(\alpha-1)}=(\alpha-1)(\rho(s)-t)^{(\alpha-2)}.
\end{equation}

\begin{equation} \label{ou2}
\nabla_t
(\rho(s)-t)^{(\alpha-1)}=-(\alpha-1)(\rho(s)-t)^{(\alpha-2)}.
\end{equation}

For the sake of the nabla fractional calculus we have the following definition

\begin{defn} \label{rising}(\cite{Adv,Boros,Grah,Spanier})

(i) For a natural number $m$, the $m$ rising (ascending) factorial of $t$ is defined by

\begin{equation}\label{rising 1}
    t^{\overline{m}}= \prod_{k=0}^{m-1}(t+k),~~~t^{\overline{0}}=1.
\end{equation}

(ii) For any real number the $\alpha$ rising function is defined by
\begin{equation}\label{alpharising}
 t^{\overline{\alpha}}=\frac{\Gamma(t+\alpha)}{\Gamma(t)},~~~t \in \mathbb{R}-~\{...,-2,-1,0\},~~0^{\mathbb{\alpha}}=0
\end{equation}

\end{defn}

Regarding the rising factorial function we observe the following:

(i) \begin{equation}\label{oper}
    \nabla (t^{\overline{\alpha}})=\alpha t^{\overline{\alpha-1}}
\end{equation}

 (ii)
 \begin{equation}\label{oper2}
    (t^{\overline{\alpha}})=(t+\alpha-1)^{(\alpha)}.
 \end{equation}

(iii)
\begin{equation}\label{oper3}
   \Delta_t (s-\rho(t))^{\overline{\alpha}}= -\alpha  (s-\rho(t))^{\overline{\alpha-1}}
\end{equation}

\textbf{Notation}:
\begin{enumerate}
\item[$(i)$] For a real $\alpha>0$, we set $n=[\alpha]+1$, where $[\alpha]$ is the greatest integer less than $\alpha$.

\item[$(ii)$] For real numbers $a$ and $b$, we denote $\mathbb{N}_a=\{a,a+1,...\}$ and $~_{b}\mathbb{N}=\{b,b-1,...\}$.

\item[$(iii)$]For $n \in \mathbb{N}$ and real $a$, we denote
$$ _{\circleddash}\Delta^n f(t)\triangleq (-1)^n\Delta^n f(t).$$
\item[$(iv)$]For $n \in \mathbb{N}$ and real $b$, we denote
                   $$ \nabla_{\circleddash}^n f(t)\triangleq (-1)^n\nabla^n f(t).$$
\end{enumerate}

\section{Definitions and essential lemmas}

\begin{defn} \label{fractional sums}
Let $\sigma(t)=t+1$ and $\rho(t)=t-1$ be the forward and backward jumping operators, respectively. Then

(i) The (delta) left fractional sum of order $\alpha>0$ (starting from $a$) is defined by:
\begin{equation}\label{dls}
    \Delta_a^{-\alpha} f(t)=\frac{1}{\Gamma(\alpha)} \sum_{s=a}^{t-\alpha}(t-\sigma(s))^{(\alpha-1)}f(s),~~t \in \mathbb{N}_{a+\alpha}.
\end{equation}

(ii) The (delta) right fractional sum of order $\alpha>0$ (ending at  $b$) is defined by:
\begin{equation}\label{drs}
   ~_{b}\Delta^{-\alpha} f(t)=\frac{1}{\Gamma(\alpha)} \sum_{s=t+\alpha}^{b}(s-\sigma(t))^{(\alpha-1)}f(s)=\frac{1}{\Gamma(\alpha)} \sum_{s=t+\alpha}^{b}(\rho(s)-t)^{(\alpha-1)}f(s),~~t \in ~_{b-\alpha}\mathbb{N}.
\end{equation}

(iii) The (nabla) left fractional sum of order $\alpha>0$ (starting from $a$) is defined by:
\begin{equation}\label{nlf}
  \nabla_a^{-\alpha} f(t)=\frac{1}{\Gamma(\alpha)} \sum_{s=a+1}^t(t-\rho(s))^{\overline{\alpha-1}}f(s),~~t \in \mathbb{N}_{a+1}.
\end{equation}

(iv)The (nabla) right fractional sum of order $\alpha>0$ (ending at $b$) is defined by:
\begin{equation}\label{nrs}
   ~_{b}\nabla^{-\alpha} f(t)=\frac{1}{\Gamma(\alpha)} \sum_{s=t}^{b-1}(s-\rho(t))^{\overline{\alpha-1}}f(s)=\frac{1}{\Gamma(\alpha)} \sum_{s=t}^{b-1}(\sigma(s)-t)^{\overline{\alpha-1}}f(s),~~t \in ~_{b-1}\mathbb{N}.
\end{equation}
\end{defn}

Regarding the delta left fractional sum we observe the following:

(i) $\Delta_a^{-\alpha}$ maps functions defined on $\mathbb{N}_a$ to
functions defined on $\mathbb{N}_{a+\alpha}$.

(ii) $u(t)=\Delta_a^{-n}f(t)$, $n \in \mathbb{N}$, satisfies the
initial value problem
\begin{equation} \label{ivpf}
\Delta^n u(t)=f(t),~~t\in N_a,~u(a+j-1)=0,~ j=1,2,...,n.
\end{equation}

(iii) The Cauchy function $\frac{(t-\sigma(s))^{(n-1)}}{(n-1)!}$
vanishes at $s=t-(n-1),...,t-1$.

\indent

Regarding the delta right fractional sum we observe the following:

(i)  $~_{b}\Delta^{-\alpha}$ maps functions defined on $_{b}\mathbb{N}$ to
functions defined on $_{b-\alpha}\mathbb{N}$.

(ii) $u(t)=~_{b}\Delta^{-n}f(t)$, $n \in \mathbb{N}$, satisfies the
initial value problem
\begin{equation} \label{ivpb}
\nabla_\ominus^n u(t)=f(t),~~t\in~ _{b}N,~u(b-j+1)=0,~ j=1,2,...,n.
\end{equation}

(iii) the Cauchy function $\frac{(\rho(s)-t)^{(n-1)}}{(n-1)!}$
vanishes at $s=t+1,t+2,...,t+(n-1)$.

\indent

Regarding the nabla left fractional sum we observe the following:

(i) $ \nabla_a^{-\alpha}$ maps functions defined on $\mathbb{N}_a$ to functions defined on $\mathbb{N}_{a}$.

(ii)$ \nabla_a^{-n}f(t)$ satisfies the n-th order discrete initial value problem

\begin{equation}\label{s1}
    \nabla^n y(t)=f(t),~~~\nabla^i y(a)=0,~~i=0,1,...,n-1
\end{equation}

(iii) The Cauchy function $\frac{(t-\rho(s))^{\overline{n-1}}}{\Gamma(n)}$ satisfies $\nabla^n y(t)=0$.

\indent

Regarding the nabla right fractional sum we observe the following:

(i) $ ~_{b}\nabla^{-\alpha}$ maps functions defined on $~_{b}\mathbb{N}$ to functions defined on $~_{b}\mathbb{N}$.

(ii)$ ~_{b}\nabla^{-n}f(t)$ satisfies the n-th order discrete initial value problem

\begin{equation}\label{s2}
    ~_{\ominus}\Delta^n y(t)=f(t),~~~ ~_{\ominus}\Delta^i y(b)=0,~~i=0,1,...,n-1.
\end{equation}
The proof can be done inductively. Namely, assuming it is true for $n$, we have
\begin{equation}\label{t1}
    ~_{\ominus}\Delta^{n+1} ~_{b}\nabla^{-(n+1)}f(t)=~_{\ominus}\Delta^{n}[-\Delta ~_{b}\nabla^{-(n+1)}f(t)].
\end{equation}

By the help of (\ref{oper3}), it follows that
\begin{equation}\label{t2}
  ~_{\ominus}\Delta^{n+1} ~_{b}\nabla^{-(n+1)}f(t)=  ~_{\ominus}\Delta^{n} ~_{b}\nabla^{-n}f(t)=f(t).
\end{equation}
The other part is clear by using the convention that $\sum_{k=s+1}^s=0$.

(iii) The Cauchy function $\frac{(s-\rho(t))^{\overline{n-1}}}{\Gamma(n)}$ satisfies $_{\ominus}\Delta^n y(t)=0$.

\indent

\begin{defn} \label{fractional differences}
(i)\cite{Miller}  The (delta) left fractional difference of order $\alpha>0$ (starting from $a$ ) is defined by:
\begin{equation}\label{dls}
    \Delta_a^{\alpha} f(t)=\Delta^n \Delta_a^{-(n-\alpha)} f(t)= \frac{\Delta^n}{\Gamma(n-\alpha)} \sum_{s=a}^{t-(n-\alpha)}(t-\sigma(s))^{(n-\alpha-1)}f(s),~~t \in \mathbb{N}_{a+(n-\alpha)}
\end{equation}

(ii) \cite{TDbyparts} The (delta) right fractional difference of order $\alpha>0$ (ending at  $b$ ) is defined by:
\begin{equation}\label{drd}
   ~_{b}\Delta^{\alpha} f(t)=  \nabla_{\circleddash}^n  ~_{b}\Delta^{-(n-\alpha)}f(t)=\frac{(-1)^n \nabla ^n}{\Gamma(n-\alpha)} \sum_{s=t+(n-\alpha)}^{b}(s-\sigma(t))^{(n-\alpha-1)}f(s),~~t \in ~_{b-(n-\alpha)}\mathbb{N}
\end{equation}

(iii) The (nabla) left fractional difference of order $\alpha>0$ (starting from $a$ ) is defined by:
\begin{equation}\label{nld}
  \nabla_a^{\alpha} f(t)=\nabla^n \nabla_a^{-(n-\alpha)}f(t)= \frac{\nabla^n}{\Gamma(n-\alpha)} \sum_{s=a+1}^t(t-\rho(s))^{\overline{n-\alpha-1}}f(s),~~t \in \mathbb{N}_{a+1}
\end{equation}

(iv) The (nabla) right fractional difference of order $\alpha>0$ (ending at $b$ ) is defined by:
\begin{equation}\label{nrd}
   ~_{b}\nabla^{\alpha} f(t)= ~_{\circleddash}\Delta^n ~_{b}\nabla^{-(n-\alpha)}f(t) =\frac{(-1)^n\Delta^n}{\Gamma(n-\alpha)} \sum_{s=t}^{b-1}(s-\rho(t))^{\overline{n-\alpha-1}}f(s),~~t \in ~ _{b-1}\mathbb{N}
\end{equation}

\end{defn}

Regarding the domains of the fractional type differences we observe:

(i) The delta left fractional difference $\Delta_a^\alpha$ maps functions defined on $\mathbb{N}_a$ to functions defined on $\mathbb{N}_{a+(n-\alpha)}$.

(ii) The delta right fractional difference $~_{b}\Delta^\alpha$ maps functions defined on $~_{b}\mathbb{N}$ to functions defined on $~_{b-(n-\alpha)}\mathbb{N}$.

(iii) The nabla left fractional difference $\nabla_a^\alpha$ maps functions defined on $\mathbb{N}_a$ to functions defined on $\mathbb{N}_{a+n}$ .

(iv)  The nabla right fractional difference $~_{b}\nabla^\alpha$ maps functions defined on $~_{b}\mathbb{N}$ to functions defined on $~_{b-n}\mathbb{N}$ .

 \begin{lem} \label{ATO} \cite{Ferd}
 For any $\alpha >0$, the following equality holds:
 $$ \Delta_a^{-\alpha} \Delta f(t)=  \Delta \Delta_a^{-\alpha}f(t)-\frac{(t-a)^{\overline{\alpha-1}}}  {\Gamma(\alpha)} f(a).$$
 \end{lem}

 \begin{lem} \label{TD} \cite{TDbyparts}
 For any $\alpha >0$, the following equality holds:
 $$~_{b} \Delta^{-\alpha} \nabla_{\circleddash} f(t)=  \nabla_{\circleddash} ~_{b} \Delta^{-\alpha}f(t)-\frac{(b-t)^{\overline{\alpha-1}}}  {\Gamma(\alpha)} f(b).$$
 \end{lem}

\begin{lem} \label{At} \cite{Gronwall}
For any $\alpha >0$, the following equality holds:
\begin{equation} \label{At1}
\nabla _{a+1}^{-\alpha} \nabla f(t)= \nabla \nabla_a^{-\alpha}f(t)-\frac{(t-a+1)^{\overline{\alpha-1}}  }{\Gamma(\alpha)}f(a)
\end{equation}
\end{lem}

 The result of Lemma \ref{At}  was obtained in \cite{Gronwall} by applying the nabla left fractional sum starting from $a$ not from $a+1$. Next will provide the version of Lemma \ref{At} by applying the definition in this article. Actually, the nabla fractional sums defined in this article and those in \cite{Gronwall} are related. For more details we refer to \cite{THFer}.

\begin{lem} \label{AtT}
For any $\alpha >0$, the following equality holds:
\begin{equation} \label{AtT1}
\nabla _a^{-\alpha} \nabla f(t)= \nabla \nabla_a^{-\alpha}f(t)-\frac{(t-a)^{\overline{\alpha-1}}  }{\Gamma(\alpha)}f(a).
\end{equation}
\end{lem}

\begin{proof}
By the help of the following by parts identity

\begin{equation}\label{AtT2}
\nabla_s [(t-s)^{\overline{\alpha-1}}f(s)]=\nabla_s (t-s)^{\overline{\alpha-1} }f(s)+ (t-\rho(s))^{\overline{\alpha-1}}  \nabla_s f(s)
\end{equation}
$$=-(\alpha-1)(t-\rho(s))^{\overline{\alpha-2}}f(s)+ (t-\rho(s))^{\overline{\alpha-1}}  \nabla_s f(s)$$
 we have
 $$\nabla_a^{-\alpha} \nabla f(t)=\frac{1}{\Gamma(\alpha)}\sum_{s=a+1}^t (t-\rho(s))^{\overline{\alpha-1}}  \nabla_s f(s)=$$
 \begin{equation}\label{AtT3}
 \frac{1}{\Gamma(\alpha)}[(t-s)^{\overline{\alpha-1} }f(s)|_a^t+(\alpha-1)\sum_{s=a+1}^t (t-\rho(s))^{\overline{\alpha-2}}f(s)]=
\end{equation}
\begin{equation}\label{AtT4}
 -\frac{(t-a)^{ \overline{\alpha-1} } } {\Gamma(\alpha)} f(a)+ \frac{1}{\Gamma(\alpha-1)}\sum_{s=a+1}^t (t-\rho(s))^{\overline{\alpha-2}}f(s)
\end{equation}

On the other hand
$$ \nabla \nabla_a^{-\alpha}f(t)=$$
\begin{equation}\label{AtT5}
 \frac{1}{\Gamma(\alpha)} \sum_{s=a+1}^t \nabla_t  (t-\rho(s))^{\overline{\alpha-1}}f(s)=\frac{1}{\Gamma(\alpha-1)}\sum_{s=a+1}^t (t-\rho(s))^{\overline{\alpha-2}}f(s)
\end{equation}

\end{proof}

\begin{rem}\label{lforany}
Let $\alpha>0$ and $n=[\alpha]+1$. Then, by the help of Lemma
 \ref{AtT} we  have
\begin{equation}\label{lforany1}
\nabla \nabla_a^\alpha f(t)=\nabla \nabla^n(\nabla_a^{-(n-\alpha)}f(t))= \nabla^n (\nabla \nabla_a^{-(n-\alpha)}f(t)).
\end{equation}
or
\begin{equation}\label{lforany11}
\nabla \nabla_a^\alpha f(t)=\nabla^n [\nabla_a^{-(n-\alpha)}\nabla
f(t)+\frac{(t-a)^{\overline{n-\alpha-1}}} {\Gamma(n-\alpha)}f(a)]
\end{equation}
Then, using the identity
\begin{equation}\label{forany2}
\nabla^n \frac{(t-a)^{\overline{n-\alpha-1}}} {\Gamma(n-\alpha)}=
\frac{(t-a)^{\overline{-\alpha-1}}} {\Gamma(-\alpha)}
\end{equation}
we infer that (\ref{AtT1}) is valid for any real $\alpha$.
\end{rem}

By the help of Lemma \ref{AtT}, Remark \ref{lforany} and the identity $\nabla (t-a)^{\overline{\alpha-1}}=(\alpha-1)(t-a)^{\overline{\alpha-2}}$, we arrive inductively at the following generalization.

 \begin{thm} \label{LNg}
 For any real number $\alpha$ and any positive integer $p$, the
following equality holds:
\begin{equation} \label{LNg1}
\nabla_a^{-\alpha}~\nabla^p f(t)=\nabla^p
\nabla_a^{-\alpha}f(t)-\sum_{k=0}^{p-1}\frac{(t-a)^{\overline{\alpha-p+k}}}
{\Gamma(\alpha+k-p+1)}\nabla^k f(a).
\end{equation}
where $f$ is defined on $\mathbb{N}_a$ and some points before $a$ .
 \end{thm}

\begin{lem} \label{RN}
For any $\alpha >0$, the following equality holds:
\begin{equation}\label{RN1}
    _{b}\nabla^{-\alpha}~ _{\circleddash}\Delta f(t)= ~ _{\circleddash}\Delta~ _{b}\nabla^{-\alpha} f(t)-\frac{(b-t)^{\overline{\alpha-1}}}  {\Gamma(\alpha)} f(b)
\end{equation}
\end{lem}

\begin{proof}
By the help of the following discrete by parts formula:
$$\Delta_s [(\rho(s)-\rho(t))^{\overline{\alpha-1}} f(s)] =$$
\begin{equation}\label{RN2}
 (\alpha-1)(s-\rho(t))^{\overline{\alpha-2}}f(s)+(s-\rho(t))^{\overline{\alpha-1}}\Delta f(s)
\end{equation}
we have
$$~_{b}\nabla^{-\alpha}~ _{a}\Delta f(t)=-\frac{1}{\Gamma(\alpha)} \sum_{s=t}^{b-1}(s-\rho(t))^{\overline{\alpha-1}}\Delta f(s)=$$

\begin{equation}\label{RN2}
  \frac{1}{\Gamma(\alpha)}[ -\sum_{s=t}^{b-1} \Delta_s ((\rho(s)-\rho(t))^{\overline{\alpha-1}} f(s) )+ (\alpha-1)\sum_{s=t}^{b-1} (s-\rho(t))^{\overline{\alpha-2}} f(s)]=
\end{equation}

\begin{equation}\label{RN3}
\frac{1}{\Gamma(\alpha-1)}\sum_{s=t}^{b-1} (s-\rho(t))^{\overline{\alpha-2}} f(s)- \frac{(b-t)^{\overline{\alpha-1}}}{\Gamma(\alpha)}  f(b).
\end{equation}

On the other hand,
$$_{\circleddash}\Delta~ _{b}\nabla^{-\alpha} f(t)=$$

\begin{equation}\label{RN3}
   - \frac{1}{\Gamma(\alpha)}\sum_{s=t}^{b-1} \Delta_t (s-\rho(t))^{\overline{\alpha-1}} f(s)=\frac{1}{\Gamma(\alpha-1)}\sum_{s=t}^{b-1} (s-\rho(t))^{\overline{\alpha-2}} f(s)
\end{equation}

where the identity

$$ \Delta_t (s-\rho(t))^{\overline{\alpha-1}}=-(\alpha-1)(s-\rho(t))^{\overline{\alpha-2}}$$
and the convention that $(0)^{\overline{\alpha-1}}=0$ are used.
\end{proof}

\begin{rem}\label{forany}
Let $\alpha>0$ and $n=[\alpha]+1$. Then, by the help of Lemma
 \ref{RN} we can have
\begin{equation}\label{forany1}
~_{a}\Delta ~_{b}\nabla^\alpha f(t)=~_{a}\Delta
~_{\circleddash}\Delta^n(~_{b}\nabla^{-(n-\alpha)}f(t))=~_{\circleddash}\Delta^n (~_{\circleddash}\Delta ~_{b}\nabla^{-(n-\alpha)}f(t))
\end{equation}
or
\begin{equation}\label{forany11}
~_{\circleddash}\Delta ~_{b}\nabla^\alpha f(t)=~_{\circleddash}\Delta^n [~_{b}\nabla^{-(n-\alpha)}~_{\circleddash}\Delta
f(t)+\frac{(b-t)^{\overline{n-\alpha-1}}} {\Gamma(n-\alpha)}f(b)]
\end{equation}
Then, using the identity
\begin{equation}\label{forany2}
~_{\circleddash}\Delta^n \frac{(b-t)^{\overline{n-\alpha-1}}} {\Gamma(n-\alpha)}=
\frac{(b-t)^{\overline{-\alpha-1}}} {\Gamma(-\alpha)}
\end{equation}
we infer that (\ref{RN1}) is valid for any real $\alpha$.
\end{rem}

By the help of Lemma \ref{RN}, Remark \ref{forany} and the identity $\Delta (b-t)^{\overline{\alpha-1}}=-(\alpha-1)(b-t)^{\overline{\alpha-2}}$, if we follow inductively we arrive at the following generalization

 \begin{thm} \label{RNg}
 For any real number $\alpha$ and any positive integer $p$, the
following equality holds:
\begin{equation} \label{RNg1}
~_{b}\nabla^{-\alpha}~_{\circleddash}\Delta^p f(t)=~_{\circleddash}\Delta^p
~_{b}\nabla^{-\alpha}f(t)-\sum_{k=0}^{p-1}\frac{(b-t)^{\overline{\alpha-p+k}}}
{\Gamma(\alpha+k-p+1)}~_{\ominus}\Delta^k f(b)
\end{equation}
where $f$ is defined on $_{b}N$ and some points after $b$.
 \end{thm}

The following theorem modifies Theorem \ref{LNg} when $f$ is only defined at $\mathbb{N}_a$.
 \begin{thm} \label{modleft}
 For any real number $\alpha$ and any positive integer $p$, the
following equality holds:
\begin{equation} \label{dLNg1}
\nabla_{a+p-1}^{-\alpha}~\nabla^p f(t)=\nabla^p
\nabla_{a+p-1}^{-\alpha}f(t)-\sum_{k=0}^{p-1}\frac{(t-(a+p-1))^{\overline{\alpha-p+k}}}
{\Gamma(\alpha+k-p+1)}\nabla^k f(a+p-1).
\end{equation}
where $f$ is defined on only $\mathbb{N}_a$ .
 \end{thm}
 The proof follows by applying   Remark \ref{lforany} inductively.

 \indent

 Similarly, in the right case we have
 \begin{thm}  \label{modright}
 For any real number $\alpha$ and any positive integer $p$, the
following equality holds:
 \begin{equation} \label{dRNg1}
~_{b-p+1}\nabla^{-\alpha}~_{\circleddash}\Delta^p f(t)=~_{\circleddash}\Delta^p
~_{b-p+1}\nabla^{-\alpha}f(t)-\sum_{k=0}^{p-1}\frac{(b-p+1-t)^{\overline{\alpha-p+k}}}
{\Gamma(\alpha+k-p+1)}~_{\ominus}\Delta^k f(b-p+1)
\end{equation}
where $f$ is defined on $_{b}N$ only.
 \end{thm}
\indent

\section{Dual identities for right fractional sums and  differences}

The dual relations for left fractional sums and differences were investigated in \cite{Nabla}. Indeed, the following two lemmas are dual relations between the delta left fractional sums (differences) and the nabla left fractional sums (differences).

\begin{lem} \label{left dual}\cite{Nabla}
Let $0\leq n-1< \alpha \leq n$ and let $y(t)$ be defined on $\mathbb{N}_a$. Then the following statements are valid.

(i)$ (\Delta_a^\alpha) y(t-\alpha)= \nabla_a^\alpha y(t)$ for $t \in \mathbb{N}_{n+a}$.

(ii) $ (\Delta_a^{-\alpha}) y(t+\alpha)= \nabla_a^{-\alpha} y(t)$ for $t \in \mathbb{N}_a$.
\end{lem}

\begin{lem} \label{left dual2}\cite{Nabla}
Let $0\leq n-1< \alpha \leq n$ and let $y(t)$ be defined on $\mathbb{N}_{\alpha-n}$. Then the following statements are valid.

(i)$ \Delta_{\alpha-n}^\alpha y(t)= (\nabla_{\alpha-n}^\alpha y)(t+\alpha)$ for $t \in \mathbb{N}_{-n}$.

(ii) $ \Delta_{\alpha-n}^{-(n-\alpha)} y(t)= (\nabla_{\alpha-n}^{-(n-\alpha)} y)(t-n+\alpha)$ for $t \in \mathbb{N}_0$.
\end{lem}
 We remind that the above two dual lemmas for left fractional sums and differences were obtained when the nabla left fractional sum was defined by

 \begin{equation}\label{remind}
  \nabla_a^{-\alpha} f(t)=\frac{1}{\Gamma(\alpha)} \sum_{s=a}^t(t-\rho(s))^{\overline{\alpha-1}}f(s),~~t \in \mathbb{N}_{a}
\end{equation}

Now, in analogous to Lemma \ref{left dual} and Lemma \ref{left dual2}, for the right fractional summations and differences we obtain

\begin{lem} \label{right dual}
Let $y(t)$ be defined on $~_{b+1}\mathbb{N}$. Then the following statements are valid.

(i)$ (~_{b}\Delta^\alpha) y(t+\alpha)= ~_{b+1}\nabla^\alpha y(t)$ for $t \in ~_{b-n}\mathbb{N}$.

(ii) $ (~_{b}\Delta^{-\alpha}) y(t-\alpha)= ~_{b+1}\nabla^{-\alpha} y(t)$ for $t \in ~_{b}\mathbb{N}$.
\end{lem}

\begin{proof}
We prove only (i). The proof of (ii) is similar and easier.
$$(~_{b}\Delta^\alpha) y(t+\alpha)=(-1)^n \nabla^n ~_{b}\Delta^{-(n-\alpha)}y(t+\alpha)=$$

\begin{equation}\label{rd1}
    \frac{(-1)^n \nabla^n}{\Gamma(n-\alpha)} \sum_{s=t+n}^{b}(s-t-1-\alpha)^{(n-\alpha-1)}y(s)= \frac{(-1)^n \Delta^n}{\Gamma(n-\alpha)} \sum_{s=t}^{b}(s-t-1-+n-\alpha)^{(n-\alpha-1)}y(s)
\end{equation}

Using the identity $t^{\overline{\alpha}}=(t+\alpha-1)^{(\alpha)}$, we arrive at
$$(~_{b}\Delta^\alpha) y(t+\alpha)=\frac{(-1)^n \Delta^n}{\Gamma(n-\alpha)} \sum_{s=t}^{b}(s-\rho(t))^{\overline{n-\alpha-1}}y(s)=$$
\begin{equation}\label{rd2}
 (-1)^n \Delta^n~ _{b+1}\nabla^{-(n-\alpha)} y(t)=  ~_{b+1}\nabla^\alpha y(t).
\end{equation}
\end{proof}

\begin{lem} \label{right dual2}
Let $0\leq n-1< \alpha \leq n$ and let $y(t)$ be defined on $~_{n-\alpha}\mathbb{N}$. Then the following statements are valid.

(i) $$~_{n-\alpha}\Delta^\alpha y(t)=~_{n-\alpha+1}\nabla^\alpha y(t-\alpha),~~~t \in ~_{n}\mathbb{N}$$

(ii)$$~_{n-\alpha}\Delta^{-(n-\alpha)} y(t)=~_{n-\alpha+1}\nabla^{-(n-\alpha)} y(t+n-\alpha),~~~t \in ~_{0}\mathbb{N}$$
\end{lem}

\begin{proof} We prove (i), the proof of (ii) is similar. By the definition of right nabla difference  we have

$$~_{n-\alpha+1}\nabla^\alpha y(t-\alpha)=~_{a}\Delta^n \frac{1}{\Gamma(n-\alpha)}\sum_{s=t-\alpha}^{n-\alpha}(s-\rho(t-\alpha))^{\overline{n-\alpha-1}}y(s)=$$

\begin{equation}\label{rig1}
  ~_{a}\Delta^n \frac{1}{\Gamma(n-\alpha)}\sum_{s=t-\alpha}^{n-\alpha}(s-\rho(t-\alpha))^{\overline{n-\alpha-1}}y(s)= \nabla_b^n \frac{1}{\Gamma(n-\alpha)}\sum_{s=t+n-\alpha}^{n-\alpha}(s-\rho(t+n-\alpha))^{\overline{n-\alpha-1}}y(s)
\end{equation}
By using (\ref{oper2}) it follows that

\begin{equation}\label{rig2}
~_{n-\alpha+1}\nabla^\alpha y(t-\alpha)=\nabla_b^n \frac{1}{\Gamma(n-\alpha)}\sum_{s=t+n-\alpha}^{n-\alpha}(s-\sigma(t))^{(n-\alpha-1)}y(s)=~_{n-\alpha}\Delta^\alpha y(t)
\end{equation}
\end{proof}
Note that the above two dual lemmas for right fractional  differences can not be obtained if we apply the definition of the delta right fractional difference introduced in \cite{Nuno} and \cite{Atmodel}.

\begin{lem}\label{power} \cite{TDbyparts}
Let $\alpha>0,~\mu>0$. Then,
\begin{equation} \label{power1}
~_{b-\mu}\Delta^{-\alpha}(b-t)^{(\mu)}=\frac{\Gamma(\mu+1)}
{\Gamma(\mu+\alpha+1)}(b-t)^{(\mu+\alpha)}
\end{equation}
\end{lem}
The following commutative property for delta right fractional sums is Theorem 9 in \cite{TDbyparts}.

\begin{thm}\label{bcomp}
Let $\alpha>0,~\mu>0$. Then, for all $t$ such that $t\equiv
b-(\mu+\alpha)~(mod~1)$, we have
\begin{equation}\label{bcomp1}
~_{b}\Delta^{-\alpha}[~_{b}\Delta^{-\mu}
f(t)]=~_{b}\Delta^{-(\mu+\alpha)}f(t)=~_{b}\Delta^{-\mu}[~_{b}\Delta^{-\alpha}
f(t)]
\end{equation}
where $f$ is defined on $_{b}N$.
\end{thm}

\begin{prop}\label{semi}
Let $f$ be a real valued function defined on $~_{b}\mathbb{N}$, and let $\alpha, \beta >0$. Then

\begin{equation}\label{semi1}
    ~_{b}\nabla^{-\alpha}[~_{b}\nabla^{-\beta}f(t)]= ~_{b}\nabla^{-(\alpha+\beta)}f(t)=  ~_{b}\nabla^{-\beta}[~_{b}\nabla^{-\alpha}f(t)]
\end{equation}
\end{prop}

\begin{proof}
The proof follows by applying Lemma \ref{right dual}(ii) and Theorem \ref{bcomp} above. Indeed,

$$ ~_{b}\nabla^{-\alpha}[~_{b}\nabla^{-\beta}f(t)]=  ~_{b}\nabla^{-\alpha} ~_{b-1}\Delta ^{-\beta}f(t-\beta)= $$

\begin{equation}\label{semi2}
  ~_{b-1}\Delta^{-\alpha} ~_{b-1}\Delta^{-\beta} f(t-(\alpha+\beta))=   ~_{b-1}\Delta^{-(\alpha+\beta)} f(t-(\alpha+\mu)) =~_{b}\nabla^{-(\alpha+\beta)}y(t)
\end{equation}
\end{proof}

The following power rule for nabla right fractional differences plays an important rule.

\begin{prop} \label{nabla right power}
Let $\alpha >0,~ \mu >-1$. Then, for $t \in ~_{b}\mathbb{N}$ , we have
\begin{equation}\label{pnr1}
~ _{b}\nabla ^{-\alpha} (b-t)^{\overline{\mu}}=\frac{\Gamma(\mu+1) } {\Gamma(\alpha+\mu+1)} (b-t)^{\overline{\alpha+\mu}}
\end{equation}

\end{prop}

\begin{proof}
By the dual formula (ii) of Lemma \ref{right dual}, we have
$$~ _{b}\nabla ^{-\alpha} (b-t)^{\overline{\mu}}= ~_{b-1}\Delta^{-\alpha} (b-r)^{\overline{\mu}}|_{r=t-\alpha}=$$
\begin{equation}\label{pnr2}
\frac{1}{\Gamma(\alpha)} \sum_{s=t}^{b-1}(s-t+\alpha-1)^{(\alpha-1)}(b-s)^{\overline{\mu}}.
\end{equation}
Then by the identity $t^{\overline{\alpha}}=(t+\alpha-1)^{(\alpha-1)}$ and using the change of variable $r=s-\mu+1$, it follows that
$$~ _{b}\nabla ^{-\alpha} (b-t)^{\overline{\mu}}=$$
\begin{equation}\label{pnr3}
  \frac{1}{\Gamma(\alpha)} \sum_{r=t-\mu+1}^{b-\mu} (r-\sigma(t-\alpha-\mu+1))^{(\alpha-1)} (b-r)^{\overline{\mu}}=(~_{b-\mu}\Delta^{-\alpha} (b-u)^{\overline{\mu}})|_{u=-\alpha-\mu+1+t}.
\end{equation}
Which  by Lemma \ref{power} leads to

$$~ _{b}\nabla ^{-\alpha} (b-t)^{\overline{\mu}}=$$
\begin{equation}\label{pnr4}
\frac{\Gamma(\mu+1) } {\Gamma(\alpha+\mu+1)}(b-t+\alpha+\mu-1)^{(\alpha+\mu)}=\frac{\Gamma(\mu+1) } {\Gamma(\alpha+\mu+1)} (b-t)^{\overline{\alpha+\mu}}
\end{equation}
\end{proof}

Similarly, for the nabla left fractional sum we can have the following power formula and exponent law.

\begin{prop} \label{power nabla left}
Let $\alpha >0,~ \mu >-1$. Then, for $t \in \mathbb{N}_a$ , we have
\begin{equation}\label{pnl1}
\nabla_a ^{-\alpha} (t-a)^{\overline{\mu}}=\frac{\Gamma(\mu+1) } {\Gamma(\alpha+\mu+1)} (t-a)^{\overline{\alpha+\mu}}
\end{equation}

\end{prop}

\begin{prop}\label{lsemi}
Let $f$ be a real valued function defined on $\mathbb{N}_a$, and let $\alpha, \beta >0$. Then

\begin{equation}\label{lsemi1}
    \nabla_a^{-\alpha}[\nabla_a^{-\beta}f(t)]= \nabla_a^{-(\alpha+\beta)}f(t)=  \nabla_a^{-\beta}[\nabla_a^{-\alpha}f(t)]
\end{equation}
\end{prop}

\begin{proof}
The proof can be achieved as in Theorem 2.1 \cite{Nabla}, by expressing the left hand side of (\ref{lsemi1}), interchanging the order of summation and using the power formula (\ref{pnl1}). Alternatively, the proof can be done by following as in the proof of Proposition \ref{semi} with the help of the dual formula for left fractional sum in Lemma \ref{left dual} after its arrangement according to our definitions.
\end{proof}

\section{Integration by parts for fractional sums and differences}

In this section we state the integration by parts formulas for nabla fractional sums and differences obtained in \cite{THFer}, then use the dual identities to obtain delta integration by part formulas.

\begin{prop} \label{summation by parts}\cite{THFer}
For $\alpha>0$, $a,b \in \mathbb{R}$, $f$ defined on $\mathbb{N}_a$ and $g$ defined on $~_{b}\mathbb{N}$, we have

\begin{equation}\label{sum1}
    \sum_{s=a+1}^{b-1}g(s) \nabla_a^{-\alpha} f(s)=\sum_{s=a+1}^{b-1}f(s)~_{b}\nabla^{-\alpha}g(s).
\end{equation}
\end{prop}
\begin{proof}
By the definition of the nabla left fractional sum we have
\begin{equation}\label{sum2}
   \sum_{s=a+1}^{b-1}g(s) \nabla_a^{-\alpha} f(s)=\frac{1}{\Gamma(\alpha)}\sum_{s=a+1}^{b-1} g(s)\sum_{r=a+1}^{s}(s-\rho(r))^{\overline{\alpha-1}}f(r).
\end{equation}
If we interchange the order of summation we reach at ( \ref{sum1}).

\end{proof}

By the help of Theorem \ref{LNg}, Proposition \ref{lsemi},
(\ref{s1}) and that $\nabla_a^{-(n-\alpha)}f(a)=0$, the authors in \cite{THFer} obtained the following left important tools which lead to a nabla integration by parts formula for fractional differences.

\begin{prop} \label{lbsdandds}\cite{THFer}
For $\alpha >0$, and $f$ defined in a suitable domain $\mathbb{N}_a$, we have

\begin{equation} \label{lbdse}
\nabla_a^{\alpha}  \nabla_a^{-\alpha}f(t)=f(t),
\end{equation}

\begin{equation} \label{lbsde1}
\nabla_a^{-\alpha} \nabla_a^{\alpha} f(t)=f(t),~~\texttt{when}~
\alpha\notin\mathbb{N},
 \end{equation}
 and
 \begin{equation} \label{lbsde2}
 \nabla_a^{-\alpha}  \nabla_a^{\alpha}f(t)=
f(t)-\sum_{k=0}^{n-1}\frac{(t-a)^{\overline{k}}}{k!} \nabla^kf(a),
\texttt{,when}~\alpha=n \in \mathbb{N}.
\end{equation}
\end{prop}

 By the help of Theorem \ref{RNg}, Proposition \ref{semi},
(\ref{s2}) and that $~_{b}\nabla^{-(n-\alpha)}f(b)=0$, the authors also in \cite{THFer} also obtained the following right important tool:

\begin{prop} \label{bsdandds}\cite{THFer}
For $\alpha >0$, and $f$ defined in a suitable domain $~_{b}\mathbb{N}$, we have

\begin{equation} \label{bdse}~_{b}\nabla^{\alpha}  ~_{b}\nabla^{-\alpha}f(t)=f(t),
\end{equation}

\begin{equation} \label{bsde1}
~_{b}\nabla^{-\alpha} ~_{b}\nabla^{\alpha} f(t)=f(t),~~\texttt{when}~
\alpha\notin\mathbb{N},
 \end{equation}
 and
 \begin{equation} \label{bsde2}
~_{b} \nabla^{-\alpha}  ~_{b}\nabla^{\alpha}f(t)=
f(t)-\sum_{k=0}^{n-1}\frac{(b-t)^{\overline{k}}}{k!}~_{\ominus}\Delta^kf(b)
\texttt{,when}~\alpha=n \in \mathbb{N}.
\end{equation}
\end{prop}

\begin{prop} \label{nabla bydiff} \cite{THFer}
Let $\alpha>0$ be non-integer and $a,b\in \mathbb{R}$ such that $a< b$ and $b\equiv
a~(mod~1)$.If $f$ is defined on $ _{b}N$ and $g$ is
defined on $N_a$, then
\begin{equation}\label{nabla bydiff1}
\sum_{s=a+1}^{b-1} f(s)\nabla_a^\alpha g(s)
=\sum_{s=a+1}^{b-1}g(s)~_{b}\nabla^\alpha f(s).
\end{equation}
\end{prop}

The proof was achieved by making use of Proposition \ref {summation by parts} and the tools Proposition \ref{lbsdandds} and Proposition \ref{bsdandds}.

Now by the above nabla integration by parts formulas and the dual identities in Lemma \ref{left dual} adjusted to our definitions and Lemma \ref{right dual} we can obtain delta integration by parts formulas.

\begin{prop} \label{delta by parts semmation}

Let $\alpha>0$, $a,b\in \mathbb{R}$ such that $a< b$ and $b\equiv
a~(mod~1)$. If $f$ is defined on $N_a$ and $g$ is defined on
$_{b}N$, then we have
\begin{equation}\label{byse}
\sum_{s=a+1}^{b -1}g(s)(\Delta_{a+1}
^{-\alpha}f)(s+\alpha)=\sum_{s=a+1}^{b-1}f(s) ~_{b-1}\Delta ^{-\alpha}g(s-\alpha).
\end{equation}

\end{prop}

\begin{prop} \label{delta by parts semmation}

Let $\alpha>0$ be non-integer and assume that $b\equiv a~(mod~1)$. If $f$ is defined on $ _{b}N$ and $g$ is
defined on $N_a$, then
\begin{equation}\label{bydiff1}
\sum_{s=a+1}^{b-1} f(s)\Delta_{a+1}^\alpha
g(s-\alpha)=\sum_{s=a+1}^{b-1}g(s)~_{b-1}\Delta^\alpha f(s+\alpha).
\end{equation}

\end{prop}

\section{Higher order fractional difference initial value problem within Riemann}

Let $\alpha> 0$  be a non-integer, $n =[\alpha]+1$ and $a(\alpha)=a+n-1$.  Consider the fractional initial difference equation

\begin{eqnarray} \label{Riem main}
                 \nabla^\alpha _{a(\alpha)-1}y(t)&=&f(t,y(t)),~~t=a(\alpha)+1, a(\alpha)+2,... \\
                 \nabla^{-(n-\alpha)}_{a(\alpha)-1}y(t)|_{t=a(\alpha)} &=& y(a(\alpha)) =c.   \nonumber
               \end{eqnarray}
Apply the sum operator $\nabla^{-\alpha}_{a(\alpha)} $ to both sides of (\ref{Riem main}) to get

\begin{equation}\label{Riem 1}
    \nabla^{-\alpha}_{a(\alpha)}\{  \nabla^{\alpha}_{a(\alpha)}y(t)+ \nabla^n \frac{(t-a(\alpha)+1)^{\overline{n-\alpha-1}}}{\Gamma(n-\alpha)}y(a(\alpha))  \}
    = \nabla^{-\alpha}_{a(\alpha)} f(t,y(t)).
\end{equation}

Applying (\ref{lbsde1}), we reach at

\begin{equation}\label{reach1}
    y(t)+ \nabla^{-\alpha}_{a(\alpha)}\nabla^n \frac{(t-a(\alpha)+1)^{\overline{n-\alpha-1}}}{\Gamma(n-\alpha)}y(a(\alpha))= \nabla^{-\alpha}_{a(\alpha)} f(t,y(t)).
\end{equation}
Now, set $g_\alpha(t)= \frac{(t-a(\alpha)+1)^{\overline{n-\alpha-1}}}{\Gamma(n-\alpha)}y(a(\alpha))$, then by Theorem \ref{modleft} with $p=n$ and $f(t)= g_\alpha(t)$ we have
\begin{equation}\label{reach2}
     \nabla^{-\alpha}_{a(\alpha)} \nabla^n g_\alpha(t)=\nabla^n \nabla^{-\alpha}_{a(\alpha)} g_\alpha(t)- \sum_{k=0}^{n-1} \frac{(t-a(\alpha))^{\overline{\alpha-n+k}}}{\Gamma(\alpha+k-n+1)} \nabla^k g_\alpha(a(\alpha))
\end{equation}

Noting that $\nabla^k g_\alpha (a(\alpha))= y(a(\alpha))$  for $k=0,1,..., n-1$, by the power formula (\ref{pnl1}) we conclude that

\begin{equation}\label{reach3}
 \nabla^{-\alpha}_{a(\alpha)} \nabla^n g_\alpha(t)= 0-\nabla^n \frac{(t-a(\alpha)+1)^{\overline{\alpha-1}} y(a(\alpha))}{\Gamma(\alpha)} - \sum_{k=0}^{n-1} \frac{(t-a(\alpha))^{\overline{\alpha-n+k}}}{\Gamma(\alpha+k-n+1)} y(a(\alpha)).
\end{equation}
Then the substitution in (\ref{reach1}) will lead to the following solution representation

\begin{equation}\label{reach final}
    y(t)=\nabla^n \frac{(t-a(\alpha)+1)^{\overline{\alpha-1}} y(a(\alpha))}{\Gamma(\alpha)}+ \sum_{k=0}^{n-1} \frac{(t-a(\alpha))^{\overline{\alpha-n+k}}}{\Gamma(\alpha+k-n+1)} c +\nabla^{-\alpha}_{a(\alpha)} f(t,y(t)).
\end{equation}
As a particular case, if $0< \alpha < 1$ then $a(\alpha)=a$ and hence

\begin{equation}\label{reach particular}
    y(t)= \frac{(t-a+1)^{\overline{\alpha-1}}}{\Gamma(\alpha)} y(a)+ \nabla^{-\alpha}_{a} f(t,y(t)),
\end{equation}
which is the result obtained in \cite{THFer}.

\indent

If $1<\alpha< 2$, then $n=2$ and $a(\alpha)=a+1$ and the initial value problem (\ref{Riem main}) becomes

\begin{eqnarray} \label{Riem main order 2}
                 \nabla^\alpha _{a}y(t)&=&f(t,y(t)),~~t=a+2, a+3,... \\
                 \nabla^{-(2-\alpha)}_{a}y(t)|_{t=a+1} &=& y(a+1) =c.   \nonumber
               \end{eqnarray}
From (\ref{reach final}), the solution is

\begin{equation}\label{reach final o2}
y(t)= \frac{(t-a)^{\overline{\alpha-3}} c}{\Gamma(\alpha-2)}+ \frac{(t-a-1)^{\overline{\alpha-2}} c}{\Gamma(\alpha-1)}+ \frac{(t-a-1)^{\overline{\alpha-1}} c}{\Gamma(\alpha)}   +\nabla^{-\alpha}_{a+1} f(t,y(t)).
\end{equation}
Combining the first two terms  then

\begin{equation}\label{reach final o3}
y(t)= \frac{(t-a)^{\overline{\alpha-2}} c}{\Gamma(\alpha-1)}+  \frac{(t-a-1)^{\overline{\alpha-1}} c}{\Gamma(\alpha)}   +\nabla^{-\alpha}_{a+1} f(t,y(t)).
\end{equation}

Again combining  we reach at

 \begin{equation}\label{reach final o2 final}
y(t)= \frac{(t-a)^{\overline{\alpha-1}} c}{\Gamma(\alpha)}+\nabla^{-\alpha}_{a+1} f(t,y(t)).
\end{equation}
If we proceed inductively we can state

\begin{thm}\label{sol rep}
Let $\alpha> 0$  be a non-integer, $n =[\alpha]+1$ and $a(\alpha)=a+n-1$.  Then the solution of the fractional initial difference equation

\begin{eqnarray} \label{Riem main n}
                 \nabla^\alpha _{a(\alpha)-1}y(t)&=&f(t,y(t)),~~t=a(\alpha)+1, a(\alpha)+2,... \\
                 \nabla^{-(n-\alpha)}_{a(\alpha)-1}y(t)|_{t=a(\alpha)} &=& y(a(\alpha)) =c,   \nonumber
              \end{eqnarray}
   is given by
\begin{equation}\label{reach final on final}
y(t)= \frac{(t-a(\alpha)+1)^{\overline{\alpha-1}} c}{\Gamma(\alpha)}+\nabla^{-\alpha}_{a(\alpha)} f(t,y(t)).
\end{equation}
\end{thm}
The surprise in Theorem \ref{sol rep} is  that the higher order Riemann fractional difference initial value problem of non-integer order needs only one initial condition which is not the case for the higher order difference equation (i.e of positive integer order). Also the solution is consisting of two terms which is not the case for fractional Cauchy problems into Riemann.
\section{The Q-operator and fractional difference equations}

If $f(s)$ is defined on $N_a\cap ~_{b}N$ and $a\equiv b~ (mod
~1)$ then $(Qf)(s)=f(a+b-s)$. The Q-operator generates a dual identity by which the left type and the right type fractional sums and differences are related. Using the change of variable $u=a+b-s$, in \cite{Th Caputo} it was shown  that
\begin{equation}\label{sum pr}
    \Delta_a^{-\alpha}Qf(t)= Q~_{b}\Delta^{-\alpha}f(t),
\end{equation}

and hence
\begin{equation}\label{pr}
    \Delta_a^\alpha Qf(t)= (Q ~_{b}\Delta^\alpha f)(t).
\end{equation}
The proof of (\ref{pr}) follows by the definition, (\ref{sum pr}) and by noting that

$$-Q\nabla f(t)=\Delta Qf(t).$$

Similarly, in the nabla case we have
\begin{equation}\label{npr sum}
   \nabla_a^{-\alpha}Qf(t)= Q~_{b}\nabla^{-\alpha}f(t),
 \end{equation}

and hence

\begin{equation}\label{rpr}
   \nabla_a^\alpha Qf(t)=( Q ~_{b}\nabla^\alpha f)(t).
\end{equation}

The proof of (\ref{rpr}) follows by the definition, (\ref{npr sum}) and that

$$-Q\Delta f(t)=\nabla Qf(t).$$

It is remarkable to mention that the Q-dual identity (\ref{pr})  can not be obtained if the definition of the delta right fractional difference introduced by Nuno R.O. Bastos et al. in \cite{Nuno} or by At{\i}c{\i} F. et al. in \cite{Atmodel} is used. Thus, the definition introduced in \cite{Th Caputo} and \cite{TDbyparts} is more convenience . Analogously, the Q-dual identity (\ref{rpr}) indicates that the nabla right Riemann  fractional differences presented in this article are also more convenient.

 It is clear from the above argument that, the Q-operator agrees with its continuous counterpart when applied to left and
right fractional Riemann Integrals Riemann derivatives. More
generally, this discrete version of the Q-operator can be used to
transform the discrete delay-type fractional functional difference dynamic
equations to advanced ones. For details in the continuous counterparts see
\cite{Thabet}.

\end{document}